\documentclass{gtart_h}

\def\ifplaintex{\expandafter\ifx\csname documentclass\endcsname\relax}

\def\gtp{{\mathsurround=0pt\it $\cal G\mskip-2mu$eometry \&\ 
$\cal T\!\!$opology $\cal P\!$ublications}}  

\def\recd{{\small Received:\qua\receiveddate\ifx\reviseddate\relax
\else\qquad Revised:\qua\reviseddate\fi\par}} 

\def\lognumber#1{\def\thelognumber{#1}}
\def\volumenumber#1{\def\thevolumenumber{#1}}
\def\volumeyear#1{\def\thevolumeyear{#1}}
\def\papernumber#1{\def\thepapernumber{#1}}
\def\pagenumbers#1#2{\def\startpage{#1}\def\finishpage{#2}}
\def\published#1{\def\publishdate{#1}}

\def\received#1{\def\receiveddate{#1}}
\def\revised#1{\def\reviseddate{#1}}
\def\accepted#1{\def\accepteddate{#1}}

\long\def\asciiabstract#1{\long\def\theasciiabstract{#1}}

\let\\\par\let\thelognumber\relax\let\thevolumenumber\relax
\let\thepapernumber\relax\let\thevolumeyear\relax\let\startpage\relax
\let\finishpage\relax\let\publishdate\relax\let\receiveddate\relax
\let\reviseddate\relax\let\accepteddate\relax\let\theasciititle\relax
\let\theasciiauthors\relax
\let\theasciiabstract\relax

\let\theasciiemail\relax

\ifplaintex
\font\logobig=cmssbx10 scaled 3836
\font\logomed=cmssbx10 scaled 2557
\else
\font\logobig=cmssbx10 scaled 4200
\font\logomed=cmssbx10 scaled 2800
\fi

\long\def\makeagttitle{   
\count0=\startpage
\agt\hfill      
\hbox to 45truept{\vbox to 0pt{\vglue -13truept{\logomed A\kern -.37em{\logobig 
T}\kern -.38em G}\vss}\hss}
\break
{\small Volume \thevolumenumber\ (\thevolumeyear)
\startpage--\finishpage\nl
Published: \publishdate}

\vglue .25truein

{\parskip=0pt\leftskip 0pt plus
1fil\def\\{\par\smallskip}{\Large\bf\thetitle}\par\medskip} \vglue
0.05truein

{\parskip=0pt\leftskip 0pt plus 1fil\def\\{\par}{\sc\theauthors}
\par\medskip}%
 
\vglue 0.03truein 

{\small\leftskip 25truept\rightskip 25truept{\bf Abstract}\stdspace\theabstract

{\bf AMS Classification}\stdspace\theprimaryclass
\ifx\thesecondaryclass\relax\else; \thesecondaryclass\fi\par
{\bf Keywords}\stdspace \thekeywords\par}\vglue 7truept

}   

\ifplaintex
\font\phead=cmsl9 scaled 950
\font\pnum=cmbx10 scaled 913
\font\pfoot=cmsl9 scaled 950
\headline{\vbox to 0pt{\vskip -4.5mm\line{\small\phead\ifnum
\count0=\startpage ISSN 1472-2739 (on-line) 1472-2747 (printed)
\hfill {\pnum\folio}\else\ifodd\count0\def\\{ }%
\ifx\theshorttitle\relax\thetitle\else\theshorttitle\fi\hfill{\pnum\folio}
\else\def\\{ and }{\pnum\folio}\hfill\ifx\theshortauthors\relax\theauthors
\else\theshortauthors\fi\fi\fi}\vss}}
\footline{\vbox to 0pt{\vglue 0mm\line{\small\pfoot\ifnum\count0=\startpage
\copyright\ \gtp\hfill\else
\agt, Volume \thevolumenumber\ (\thevolumeyear)\hfill\fi}\vss}}
\else
\font=cmsl9 scaled 1050
\font\lnum=cmbx10 
\font=cmsl9 scaled 1050
\makeatletter
\def\@oddhead{{\small\lhead\ifnum\count0=\startpage ISSN 1472-2739 
(on-line) 1472-2747 (printed)\hfill {\lnum\number\count0}\else\ifodd\count0
\def\\{ }\ifx\theshorttitle\relax \thetitle \else\theshorttitle\fi\hfill
{\lnum\number\count0}\else\def\\{ and }{\lnum\number\count0}
\hfill\ifx\theshortauthors\relax 
\theauthors\else\theshortauthors\fi\fi\fi}}\def\@evenhead{\@oddhead}
\def\@oddfoot{\small\lfoot\ifnum\count0=\startpage\copyright\ \gtp\hfill\else
\agt, Volume \thevolumenumber\ (\thevolumeyear)\hfill\fi}
\def\@evenfoot{\@oddfoot}
\makeatother
\fi
\let\maketitlepage\makeagttitle
\let\makeshorttitle\maketitlepage
\let\maketitle\maketitlepage

\newwrite\gtoutfile
\long\gdef\makeheadfile{  
{\def\\{, }\def\s{ }
\immediate\openout\gtoutfile head.xxx
\immediate\write\gtoutfile{Proxy-for: \ifx\theasciiauthors\relax
\theauthors\else\theasciiauthors\fi\s<\ifx\theasciiemail\relax\theemail\else\theasciiemail\fi>}
\immediate\write\gtoutfile{\noexpand\\}
\immediate\write\gtoutfile{Authors: \ifx\theasciiauthors\relax
\theauthors\else\theasciiauthors\fi}
{\def\\{ }\immediate\write\gtoutfile{Title: \ifx\theasciititle\relax
\thetitle\else\theasciititle\fi}}
\immediate\write\gtoutfile{Subj-class: GT or SG, GR etc}
\immediate\write\gtoutfile{MSC-class: \theprimaryclass\ifx\thesecondaryclass\relax\else, \thesecondaryclass\fi}
\immediate\write\gtoutfile{Journal-ref: Algebr. Geom. Topol. \thevolumenumber\s
(\thevolumeyear) \startpage-\finishpage}
\immediate\write\gtoutfile{Comments: Published by Algebraic and
Geometric Topology at}
\immediate\write\gtoutfile{\s\s\s  http://www.maths.warwick.ac.uk/agt/AGTVol\thevolumenumber/agt-\thevolumenumber-\thepapernumber.abs.html}
\immediate\write\gtoutfile{\noexpand\\}
\immediate\write\gtoutfile{}
\ifx\theasciiabstract\relax
\immediate\write\gtoutfile{\theabstract}\else
\immediate\write\gtoutfile{\theasciiabstract}\fi
\immediate\write\gtoutfile{}
\immediate\write\gtoutfile{\noexpand\\}
\immediate\write\gtoutfile{}
\immediate\closeout\gtoutfile}}  

\def\maketitlepage{\makeagttitle\makeheadfile}
\let\makeshorttitle\maketitlepage
\let\maketitle\maketitlepage

\lognumber{8}
\volumenumber{4}
\volumeyear{2004}
\papernumber{8}
\published{15 March 2004}
\pagenumbers{121}{131}
\received{19 December 2002}
\revised{5 March 2004}
\accepted{15 March 2004}

\usepackage{amsmath, amssymb}
\usepackage{amscd}

\theoremstyle{plain}
\newtheorem{theorem}{Theorem}[section]
\newtheorem{lemma}[theorem]{Lemma}

\newtheorem{proposition}[theorem]{Proposition}

\newtheorem{corollary}[theorem]{Corollary}

\theoremstyle{definition}

\newtheorem*{theorem*}{Theorem}

\theoremstyle{remark}
\newtheorem*{remark}{Remark}

\newtheorem*{example}{Examples}

\DeclareMathOperator{\Hom}{Hom}
\DeclareMathOperator{\res}{res}

\DeclareMathOperator{\Res}{Res}

\DeclareMathOperator{\Mod}{Mod}
\DeclareMathOperator{\modu}{mod}

\DeclareMathOperator{\CS}{CS}
\DeclareMathOperator{\proj}{proj}

\DeclareMathOperator{\id}{id}
\DeclareMathOperator{\add}{add}

\DeclareMathOperator{\ES}{ES}
\DeclareMathOperator{\stalks}{stalks}

\begin{document}

\title{Smith Theory for algebraic varieties}

\author{Peter Symonds}
\address{Department of Mathematics, UMIST\\PO Box 88, 
Manchester M60 1QD, UK}
\email{Peter.Symonds@umist.ac.uk}

\begin{abstract} We show how an approach to Smith Theory about group actions on CW--complexes using Bredon cohomology can be adapted to work for algebraic varieties.
\end{abstract}

\asciiabstract{We show how an approach to Smith Theory about group 
actions on CW-complexes using Bredon cohomology can be adapted to
work for algebraic varieties.}

\maketitle

\begin{section}{Introduction}
\label{intro}

Peter May described in \cite{may} a version of Smith Theory based on
Bredon cohomology, so it really applies to any complex of projective
coefficient systems rather than just a topological space. Later Jeremy
Rickard in \cite{rickard} showed how to associate a complex of
$p$-permutation modules to a group action on a variety in such a way
that the cohomology of this complex is the \'etale cohomology of the
variety. We show how to generalize this to obtain a complex of
projective coefficient systems. Thus Smith Theory becomes available
for algebraic varieties, even over fields of finite
characteristic. Our framework is also sufficient to apply to varieties
methods of Borel, Swan and others based on equivariant cohomology,
although we do not set out the details here.

\end{section}

\begin{section}{Coefficient systems}
\label{coeff}

A coefficient system $L$ on a group $G$ over a ring $R$ is a functor
from the right orbit category of $G$ to $R$-modules. In more concrete
terms, it consists of a collection of $R$-modules $L(H)$, one for each
subgroup $H \le G$ together with $R$-linear restriction maps $\res
^H_K \co L(H) \rightarrow L(K)$ for each $K \le H \le G$ and
conjugation maps $c_{g,H}\co L(H) \rightarrow L({}^gH)$ for each $g
\in G$ and $H \le G$.

These must satisfy the identities:
\newpage
\begin{enumerate}
\item
$\res ^H_H = \id, \quad H \le G$;
\item
$\res ^K_J \res ^H_K = \res ^H_J, \quad J \le K \le H \le G$;
\item
$c_{g_1,{}^{g_2}H} c_{g_2,H} = c_{g_1 g_2,H}, \quad H \le G, \quad g_1,g_2 \in G$;
\item
$\res ^{{}^gH}_{{}^gK} c_{g,H} = c_{g,K} \res ^H_K, \quad K \le H \le G,\quad g \in G$;
\item
$c_{h,H} = \id, \quad H \le G, \quad h \in H$.
\end{enumerate}

In particular, the conjugation maps make $L(H)$ into a left $RN_G(H)/H$-module.

A morphism $f\co L \rightarrow M$ is a collection of $R$-linear maps $f(H) \co L(H) \rightarrow M(H)$ which commute with the $\res$ and $c$.

The coefficient systems on $G$ over $R$ form an abelian category, which we denote by $\CS _R(G)$. If $H \leq G$ there is a forgetful map $\Res  ^G_H \co \CS_R(G) \rightarrow \CS_R(H)$

\begin{example}\begin{enumerate}
\item
The constant coefficient system $\bar R$, which is just $R$ on each evaluation and all the maps are the identity.
\item
The fixed point coefficient system $V^?$, where $V$ is a left $RG$-module and the notation indicates that the evaluation on $H \le G$ is the fixed point submodule $V^H$. Restriction is inclusion and conjugation is multiplication by $g \in G$. These have the important property that $\Hom _{\CS_R(G)}(L,V^?) \cong \Hom _{RG}(L(1),V)$.
\item
A variation on $V^?$ is $V_0$, which takes the value $V$ on $1$ and $0$ elsewhere.
\item
The systems $R[X^?]$, where $X$ is a left $G$-set and the evaluation at $H$ is the free $R$-module on the fixed point set $X^H$. 

The particular cases $R[G/H^?]$ have the important property $$\Hom _{\CS _R(G)}(R[G/H^?],L) \cong L(H).$$ It follows that they are projective and that they provide enough projectives. Thus every projective is a summand of a sum of these.
\end{enumerate}
\end{example}

For more information on coefficient systems see \cite{symonds}.

Given a set of coefficient systems $I$ it is convenient to define $\add (I)$ to be the full subcategory of $\CS_R(G)$ in which the objects are isomorphic to a summand of a coefficient system of the form $L_1 \oplus \ldots \oplus L_n, \quad L_i \in I$.

Thus the subcategory $\proj (\CS_R(G))$ of finitely generated projective coefficient systems is the same as $\add (\{ R[G/H^?]: H \le G \})$.

If $X$ is a $G$--CW--complex then there is a complex of coefficient systems $C[X^?]$ associated to it, in which $C_n[X]= R[(X_n)^?]$ where $X_n$ is the $G$-set of $n$-cells in $X$ and the boundary morphisms are defined in the usual way.

The Bredon cohomology of $X$ with coefficients in a coefficient system $L$, as defined in \cite{bredon}, is $H^*_G(X,L)=H^* (\Hom _{\CS _R(G)}(C[X^?],L))$.

\begin{example} \begin{enumerate}
\item
 $H^*_G(X,(RG)^?) \cong H^*(X,R)$, the usual CW--cohomology,
 \item
 $H^*_G(X, R_G) \cong H^*(X^G,R)$, where $R_G$ takes the value R on $G$ and 0 elsewhere.
 \item
 $H^*_G(X, \bar R) \cong H^*(X/G,R)$,
 \item
 More generally we can regard $H^*_?(X, R)$ as a coefficient system itself under the natural restriction and conjugation maps, and then we have $H^*_?(X,\bar R) \cong H^*(X/?,R)$.
\end{enumerate}
\end{example}

The dual concept to that of a coefficient system we term an efficient system, in which the restriction maps go in the opposite direction. $E(H)$ is now a right $N_G(H)$-module, although we could remedy this by taking the contragredient instead of the dual. The category of efficient systems for $G$ over $R$ is denoted by $\ES_R(G)$.

If $R$ is self-injective then applying $\Hom _R(-,R)$ provides a duality between the subcategories taking values in finitely generated modules.

The dual of $R[G/H^?]$ is denoted by $R[G/H^?]^*$. The evaluation on $K \le G$ can be thought of as the functions on the fixed point set $(G/H)^K$ and the restriction maps just restrict the functions. If $R$ is self injective then $R[G/H^?]^*$ is injective.

$\CS _R(G)$ can also be viewed as the category of modules over an $R$-algebra $C _R(G)$ of finite rank over $R$, (cf.\ \cite{bouc}). Similarly $\ES _R(G)$ is equivalent to the category of modules over another $R$-algebra $E _R(G)$.

\end{section}

\begin{section}{Varieties}
\label{var}

From now on $k$ is an algebraically closed field and in this section $X$ is a separated scheme of finite type over $k$.

Let $A$ be a torsion Artin algebra and let $\mathcal F$ be a constructible sheaf of $A$-modules over $X$. Let $\stalks ( \mathcal F)$ denote the set of stalks of $\mathcal F$ at the $k$-rational points. This contains only a finite number of isomorphism classes of $A$-modules.

Recall that $R\Gamma _c(X,\mathcal F)$ is a complex of $A$-modules, natural as an object of the derived category $D (A- \Mod )$,  whose homology is \'etale cohomology with compact supports $H^*_c(X,\mathcal F)$.

Our main tool will be the following result from \cite{rickard}:

\begin{theorem}[Rickard]
\label{rick}
There is a complex of modules in $\add ( \stalks ( \mathcal F))$ of finite type, which we denote by $\Omega _c (X , \mathcal F)$. It is well defined up to homotopy equivalence. It has the following properties:
\begin{enumerate}
\item
$\Omega _c (X, \mathcal F)$ is isomorphic to $R \Gamma _c (X, \mathcal F)$ in $D (A- \Mod )$;
\item
$\mathcal F \mapsto \Omega _c (X, \mathcal F)$ is a functor from constructible sheaves of $A$-modules over $X$ to $K^b ( A - \modu )$;
\item
If $f\co Y \rightarrow X$ is a finite morphism of separated schemes of finite type over $k$ then there is an induced map $\Omega _c (X,  \mathcal F) \rightarrow  \Omega _c (Y, f^* \mathcal F)$;
\item
If $B$ is also a torsion Artin algebra and $L$ is a functor $\add (\stalks (\mathcal F)) \rightarrow B -\modu$ then $L \Omega _c (X, \mathcal F) \cong \Omega _c (X, \tilde {L}\mathcal F)$, where $\tilde {L} \mathcal F$ denotes the sheafification of the presheaf $L \mathcal F$.
\end{enumerate}
\end{theorem}

We will apply this in the case that $R = \mathbb Z / \ell ^n$ and $A= E _R(G)$.

We suppose that a finite group $G$ acts on $X$ with quotient variety $Y = X/G$ and
projection map $\rho \co X \rightarrow Y$. We let $\mathcal F = \mathcal F _X $ be the
sheafification of the presheaf that sends a Zariski open set $U \subseteq Y$ in the Zariski
topology to $R[( \pi _0 ( \rho ^{-1}U))^?]^*$, where $\pi _0 (\rho ^{-1}U)$ is the $G$-set
of components of $\rho ^{-1}U$. (This extends to the \'etale site on $X$ by evaluating on the image of an \'etale map $U \rightarrow X$.) Then $ \stalks (\mathcal F)$ consists of injective modules.

Theorem~\ref{rick} produces a complex of injective efficient systems of finite type $\Omega _c(Y,\mathcal F)$. These complexes for different $n$ can be pieced together in such a way that we can take the inverse limit and obtain a complex of finite type of efficient systems in $\add ( \{ \hat{ \mathbb Z _ \ell }[G/H^?]^*: H \text{\, the stabilizer of a $k$-rational point} \})$ as in \cite{rickard}. The dual of this by $\Hom _{\hat{\mathbb Z _ \ell}}(-,\hat{\mathbb Z _\ell})$ is the complex that we will denote by $C[X^?]$.

\begin{theorem}
\label{prop}
For any $H \leq G$, $C[X^?](H)^* \cong R \Gamma _c (X^H, \hat{\mathbb Z _\ell})$.
\end{theorem}

In other words $C[X^?](H)$ is a complex whose dual has cohomology $H^*_c(X^H, \hat{\mathbb Z _\ell})$. We can therefore think of it as the analogue of the complex $C[X^?]$ for the Bredon cohomology of a $G$--CW--complex.

Since $C[X^H](1) \cong C[X^?](H)$ our notation is justified and, after the proof is complete, we will write $C[X^H]$ instead of $C[X^?](H)$.

We will prove theorem~\ref{prop} as a corollary of some more general results.

Notice that $C[X^?]$ is natural with respect to group homomorphisms $f \co H \rightarrow G$
for which the kernel acts trivially on $X$.

Let $\mathcal A$ be a set of subgroups of $G$ closed under supergroups and conjugation. Let
$S_{\mathcal A} X = \bigcup _{J \in \mathcal A} X^J$.

Define $L_{\mathcal A}\co \CS _R(G) \rightarrow \CS_R(G)$ by taking $L_{\mathcal A}C$ to be
the smallest subsystem of $C$ that is equal to $C(H)$ for all $H \in \mathcal A$. Then
\begin{align*}
L_{\mathcal A}R[G/H^?] &= \begin{cases} R[G/H^?] \quad H \in \mathcal A \\ 0  \quad\quad \text{otherwise} \end{cases}\\
&=R[(S_{\mathcal A}G/H)^?].
\end{align*}
So $L_{\mathcal A}$ induces a functor $L_{\mathcal A}\co \add (\{ R[G/H^?]; H \leq G \}) \rightarrow \add ( \{ R[G/H^?] ; H \in \mathcal A \})$.

\begin{proposition}
\label{loc}
$L_{\mathcal A}C[X^?]$ is homotopy equivalent to $C[ (S_{\mathcal A}X)^?]$.
\end{proposition}

\begin{proof} It is sufficient to prove the analogous statement for $R=\mathbb Z / \ell ^n$.

By \ref{rick}, $L_{\mathcal A}C[X^?]^* \cong \Omega _c(Y,\tilde L_{\mathcal A} \mathcal F_1)$, where $\mathcal F_1$ is the sheafification of $\mathcal F_1'\co U \mapsto \Gamma ((\pi _0(\rho ^{-1} U))^?,R)$. So $\tilde L_{\mathcal A} \mathcal F_1$ is the sheafification of $U \mapsto \Gamma ((S_{\mathcal A}\pi _0 (\rho ^{-1}U))^?,R)$.

Now $C[ (S_{\mathcal A}X)^?]^* \cong \Omega _c( (S_{\mathcal A}X)/G, \mathcal F_{S_{\mathcal A}X}) \cong \Omega _c (Y, \mathcal F_2)$, where $\mathcal F_2$ is the sheafification of $\mathcal F_2'\co U \mapsto \Gamma ((\pi _0(\rho ^{-1} U \cap S_{\mathcal A}X))^?,R)$.

Inclusion of fixed points gives a map $\tilde L_{\mathcal A} \mathcal F_1' \rightarrow \mathcal F_2'$, which induces an isomorphism on the stalks and hence an isomorphism of sheaves.
\end{proof}

\begin{lemma}
\label{one}
$C[X^?](1)^* \cong R \Gamma _c(X,R)$
\end{lemma}

\begin{proof}
Again it is enough to work over $\mathbb Z / \ell ^n$.

Considering the functor ``evaluate at $1$'', we find $\Omega _c(Y, \mathcal F)(1) \cong \Omega _c (Y, \widetilde{\mathcal F (1)})$, where $\widetilde{\mathcal F (1)}$ is the sheafification of $U \mapsto \Gamma (\pi _0(\rho ^{-1}U),R)$, which is just $\rho _*R$.

Finally $\Omega _c(Y,\rho _*R) \cong R \Gamma _c(Y,\rho_*R) \cong R \Gamma _c (X,R)$.
\end{proof}

\begin{remark}
The above lemma shows that $C[X]=C[X^?](1)$ is the dual of Rick\-ard's complex of $\ell$-permutation modules.
\end{remark}

\begin{proof} of \ref{prop}. We can restrict to $N_G(H)$ if necessary, by naturality under inclusions, so we may assume that $H$ is normal in $G$. Let $\mathcal H$ be the set of subgroups of $G$ containing $H$ and apply \ref{loc} to obtain $C[X^H](H) \cong L_{\mathcal H}C[X^?](H) = C[X^?](H)$.

Notice that $C[X^H](H) \cong C[X^H](1)$ by naturality under $G \rightarrow G/H$ and apply \ref{one}.
\end{proof}

A similar method will prove the following result (see \cite{rickard}). We abbreviate $H_0(G,M)$ to $M_G$.

\begin{lemma}
\label{quot}
$C[X]_G \cong C[X/G]$.
\end{lemma}

Let $\mathcal S^1$ denote the set of non-trivial subgroups and write $S=S_{\mathcal S^1}$ and $L=L_{\mathcal S^1}$.

\begin{lemma}
\label{tri}
$(C[X^?]/LC[X^?])^* \cong R\Gamma _c(X \smallsetminus SX,R)_0$.
\end{lemma}

\begin{proof}
Both sides are zero on non-trivial subgroups, so we only need to check at the trivial group.

The inclusion map $C[SX] \rightarrow C[X]$ is equivalent to $LC[X] \rightarrow C[X]$. It is also dual to $R \Gamma _c(X,R) \rightarrow R \Gamma _c(SX,R)$. Thus the triangles $LC[X] \rightarrow C[X] \rightarrow C[X]/LC[X] $ and $R\Gamma _c(X \smallsetminus SX,R) \rightarrow R \Gamma _c(X,R) \rightarrow R \Gamma _c(SX,R)$ are dual.
\end{proof}

\end{section}

\begin{section}{Smith Theory}
\label{smith}

Various results are known collectively as Smith Theory (see \cite{bredon2}, for example), but the prototype is the theorem that if a $p$--group $P$ acts on a finite dimensional CW--complex which has the mod--$p$ cohomology of a point then the fixed point subcomplex also has the mod--$p$ cohomology of a point.  Once the case of $P$ of order $p$ is proved this follows by induction on the order of $P$.

From now on we will take $R$ to be $\mathbb F_p$. We allow $X$ to be either a CW--complex, in which case our results are well known, or a separated scheme of finite type over an algebraically closed field $k$. In the latter case the $\ell$ in the previous section becomes $p$ and as a consequence we will need the characteristic of $k$ not to be equal to $p$ in order to be able to use the \'etale cohomology.

As before, we define $H^*_G(X,L)=H^* (\Hom _{\CS _R(G)}(C[X^?],L))$.

\begin{lemma}
\label{ident}
We have the following identifications:
\begin{align*}
H^*_G(X,(RG)^?) &\cong H^*_c(X,R), \\
H^*_G(X,(RG^?)/(RG)_0) &\cong H^*_c(SX,R) \\
H^*_G(X,R_0) &\cong H^*_c((X \smallsetminus SX)/G,R).
\end{align*}
\end{lemma}

The analogous result for $G$--CW--complexes is well known.

\begin{proof} By the adjointness property of $(RG)^?$,
\begin{align*}
\Hom _{\CS _R(G)}(C[X^?],(RG)^?) & \cong \Hom _{RG} (C[X],RG) \\
& \cong \Hom _R(C[X],R) \\
& \cong R \Gamma _c (X,R).
\end{align*}

Notice that $C[SX] \cong LC[X]$, by \ref{loc}. Because $LC[X]$ is in $\add( \{ R[G/H^?];$\break $H \ne 1 \})$ we obtain
\begin{align*}
\Hom _{\CS _R(G)}(C[X^?],(RG)^?/(RG)_0) & \cong \Hom _{\CS _R(G)}(LC[X^?],(RG)^?) \\
& \cong \Hom _{RG}(LC[X],RG) \\
& \cong \Hom _{RG}(C[SX],RG) \\
& \cong \Hom _R(C[SX],R) \\
& \cong R \Gamma _c (SX,R).
\end{align*}

There are no non-zero homomorphisms from $\add( \{ R[G/H^?]; H \ne 1 \} )$ to $R_0$. Also $C[X \smallsetminus SX]$ is in $\add ( \{ R[G^?] \})$, so vanishes off the trivial group. We find that
\begin{align*}
\Hom _{\CS _R(G)}(C[X^?],R_0) & \cong \Hom _{\CS _R(G)}(C[X^?]/LC[X^?], R_0) \\
& \cong \Hom _{\CS _R(G)}(C[(X \smallsetminus SX)^?], R_0) \\
& \cong \Hom _{RG}(C[X \smallsetminus SX],R) \\
& \cong R \Gamma _c((X \smallsetminus SX)/G,R),
\end{align*}
by \ref{quot}.
\end{proof}

May's approach to Smith Theory considers the Bredon cohomology groups in the lemma above and uses various long exact sequences associated to a short exact sequence of coefficient systems.

Let $I$ denote the augmentation ideal of $RG$. Notice that if $G$ is a $p$--group, which we will denote by $P$, then $(RP)^?/I_0 \cong (RP)^?/(RP)_0 \oplus R_0$ and the composition factors of $I_0$ are all $R_0$.

Let \begin{align*}
a_q &= \dim H^q_G(X,R_0) = \dim H^q_c((X \smallsetminus SX)/G,R), \\
b_q &= \dim H^q_G(X,(RG)^?) = \dim H^q_c(X,R), \\
c_q &= \dim H^q_G(X,(RG^?)/(RG)_0) = \dim H^q_c(SX,R).
\end{align*}

May proves the following result in \cite{may} for $P$--CW--complexes but, since the proof uses only manipulations with Bredon cohomology and the identifications in \ref{ident}, it is valid for separated schemes of finite type too.

\begin{theorem}[Floyd, May] 
\label{may}
The following inequality holds for any $q \ge0$ and $r \ge 0$:
$$a_q + \sum ^r_{i=0} (|P|-1)^ic_{q+i} \le \sum ^r_{i=0} (|P|-1)^i b_{q+i} + (|P| -1)^{r+1}a_{q+r+1}.$$
In particular, if $a_i=0$ for $i$ sufficiently large,
$$a_q + \sum _{i\ge 0} (|P|-1)^ic_{q+i} \le \sum _{i \ge 0} (|P|-1)^ib_{q+i}.$$
Moreover, if $a_i, b_i, c_i =0$ for $i$ sufficiently large then
$$\chi _c (X) = \chi _c(SX) + |P| \chi_c((X \smallsetminus SX)/P).$$

If $P$ is cyclic of order $p$, and $r$ is even if $p \ne 2$, then we can remove the factors $(|P|-1)$, i.e.\
$$a_q + \sum ^r_{i=0} c_{q+i} \le \sum ^r_{i=0} b_{q+i} +  a_{q+r+1}.$$
\end{theorem}

\begin{remark}
\begin{enumerate}
\item If $X$ is a CW--complex then we can use ordinary cohomology instead of compactly supported cohomology provided that we also
replace $(X \smallsetminus SX)/G$ by $(X/SX)/G$ and take its reduced cohomology.
\item Notice that the last line includes Illusie's result \cite{illusie} for varieties that
if $P$ acts freely on $X$ then $|P|$ divides $\chi _c(X)$. In fact, in this case, $C[X]$ is
a complex of projective $RP$-modules and, since $P$ is a $p$--group, the modules are free.
\item In the topological case, if we take $X$ to be $EP$ (the universal cover of the
classifying space) and $q=0$ then we recover the well-known result that the $H^i(P, \mathbb
F_p)$ are non-zero in every degree.
\end{enumerate}
\end{remark}

Recall that $$H^i_c(\mathbb A^n(k), \mathbb F_p) = \begin{cases} \mathbb F_p, \quad i=n \\
0, \quad \text{ otherwise} \end{cases}$$ provided that $p$ is not the characteristic of
$k$.

\begin{corollary}
\label{aff}
Suppose that $X$ has the cohomology of an affine space $\mathbb A^n$ and also that if $X$
is a CW--complex then it is finite-dimensional. Then $X^P$ has the cohomology
of some affine space $\mathbb A^m$ for some $m$, with $n-m$ even if $|P| \ne 2$.
\end{corollary}

\begin{remark} \begin{enumerate}
\item
By taking $n=0$ this includes the case that when $X$ is mod-$p$ acyclic then $X^P$ must also be
mod-$p$ acyclic.
\item When $X$ is compact a similar argument shows that if $X$ is a mod-$p$ homology sphere then so is $X^P$.
\end{enumerate}
\end{remark}

\begin{proof}
By induction on $P$ we can reduce to the case when $P$ is cyclic of order $p$. For $P$ must have a normal subgroup $Q$ of index $p$, and by induction $X^Q$ has the cohomology of an affine space. But $X^P=(X^Q)^P$.

From the last line in \ref{may} with $r$ large it follows that $\sum _{i\geq 0}c_i \leq 1$. The sum can not be $0$ by the Euler characteristic formula.
\end{proof}

We now present a more conceptual approach to these results which shows how coefficient
systems can provide a very flexible tool. It is based on the following lemma:

\begin{lemma}
\label{split}
Any monomorphism between two projective coefficient systems in $\CS_R(P)$ is split.
\end{lemma}

\begin{proof}
Consider a map $R[P/U^?] \rightarrow R[P/V^?]$. It must be zero unless $U$ is conjugate to
a subgroup of $V$. But then it can only be a monomorphism if $|U| \ge |V|$ so in fact $U$
is conjugate to $V$ and the map is an isomorphism.

Now any projective $F$ is of the form $F \cong \bigoplus _{j \in J} F_j$, where each
$F_j$ is an indecomposable projective, so isomorphic to some $R[P/V^?]$. So suppose that we have a monomorphism $f\co R[P/U^?] \rightarrow \bigoplus _{j \in J} F_j$. The socle of
$R[P/U^?]$ is just the sub-system generated by $\sum _{g \in P/U} gP$ in degree $0$. One of
the components of $f$, say $f_j\co R[P/U^?] \rightarrow F_j$ must be non-zero on the socle,
hence a monomorphism and so an isomorphism. The splitting is now projection onto $F_j$
followed by $(f_j)^{-1}$.

Now consider the case $f\co\bigoplus _{i \in I} E_i \rightarrow F$. If $I$ is finite, say $I = \{1, \ldots ,n \}$, then we have a proof by induction on $n$. We have shown that $F \cong E_1 \oplus F/f(E_1)$, and there is an injection $f'\co\bigoplus _{i \in I \setminus \{1 \}} E_i \rightarrow F/f(E_1)$. The latter splits by the induction hypothesis.

The case of finite $I$ is enough for us to deduce that, for any $I$, the map $f$ is pure. But, for modules over an Artin algebra, any projective module is pure injective (because it is a summand of a free module and the free module of rank 1 is $\Sigma$-pure-injective by condition (iii) of theorem 8.1 in \cite{jl}), so $f$ is split.
 \end{proof}

\begin{corollary}
\label{splitexact}
If $C$ is a complex of projectives in $\CS_R(P)$ that is bounded above and such that $C(1)$ is exact then $C$ is split exact.
\end{corollary}

\begin{proof} Evaluation at $1$ detects monomorphisms between projectives, so there is an easy argument based on \ref{split} and induction on the number of boundary maps that can be split, starting from the left.
\end{proof}

\begin{corollary} Let $f\co C \rightarrow D$ be a map between two bounded complexes of projective coefficient systems in $\CS _R(P)$. If $f(1)\co C(1) \rightarrow D(1)$ is a quasi-isomorphism, i.e.\ induces an isomorphism in homology, then $f$ is a homotopy equivalence.
\end{corollary}

\begin{proof}
Apply \ref{splitexact} to the cone of $f$, to deduce that $f$ is a quasi-isomorphism. Since the complexes consist of projectives, $f$ must be a homotopy equivalence.
\end{proof}

\begin{corollary}
Let $f\co X \rightarrow Y$ be a finite morphism of separated schemes that induces an isomorphism on \'etale cohomology with coefficients in $R$ . Suppose that $P$ acts on both $X$ and $Y$ and that $f$ is equivariant. Then the induced morphism $f^P\co X^P \rightarrow Y^P$ also induces an isomorphism on cohomology.
\end{corollary}

\begin{remark}
It is not sufficient to consider complexes of $p$-permutation modules. For example, if we let $C_2$ denote the cyclic group of order 2 take $R = \mathbb F_2$ then there is a short exact sequence $R \rightarrow RC_2 \rightarrow R$. But this is not split.
\end{remark}

\begin{remark}
The methods of equivariant cohomology of Borel \cite{borel} can also be applied to varieties. They all depend on analyzing the triangle $C[S_{\mathcal A}X] \rightarrow C[X] \rightarrow C[X]/L_{\mathcal A}C[X]$ of $RG$-modules. The proofs in \cite{brown} and \cite{tomdieck} carry over.
\end{remark}

\end{section}

\Addresses\recd
\end{document}